\documentclass[12pt]{article}
\usepackage{latexsym}

\def\bs{\begin{subequations}}
\def\es{\end{subequations}}

\catcode`\@=11

\newtoks\@stequation
\def\subequations{\refstepcounter{equation}
  \edef\@savedequation{\the\c@equation}%
  \@stequation=\expandafter{\theequation}%   %only want \theequation
  \edef\@savedtheequation{\the\@stequation}% % expanded once
  \edef\oldtheequation{\theequation}%
  \setcounter{equation}{0}%
  \def\theequation{\oldtheequation\alph{equation}}}

\def\endsubequations{\setcounter{equation}{\@savedequation}%
  \@stequation=\expandafter{\@savedtheequation}%
  \edef\theequation{\the\@stequation}\global\@ignoretrue}

\makeatletter%  allow access to internal LaTeX commands
        \renewcommand{\theequation}{\thesection.\arabic{equation}}%
        \@addtoreset{equation}{section}%
\makeatother%  turn off access to internal LaTeX commands

\renewcommand{\thefootnote}{\fnsymbol{footnote}}

\begin{document}

\begin{titlepage}

July 7, 2008, added Appendix A

\begin{center}        \hfill   \\
            \hfill     \\
                                \hfill   \\

\vskip .25in

{\large \bf Analytic Functions of a General Matrix Variable \\}

\vskip 0.3in

Charles Schwartz\footnote{E-mail: schwartz@physics.berkeley.edu}

\vskip 0.15in

{\em Department of Physics,
     University of California\\
     Berkeley, California 94720}
        
\end{center}

\vskip .3in

\vfill

\begin{abstract}

Recent innovations on the  differential calculus for functions 
of non-commuting variables, begun for a quaternionic variable,
are now extended to the case of a general matrix over 
the complex numbers. The expansion of F(X+Delta) is given to 
first order in Delta for general matrix variables X and Delta that 
do not commute with each other. 

\end{abstract}

\vfill

\end{titlepage}

\renewcommand{\thefootnote}{\arabic{footnote}}
\setcounter{footnote}{0}
\renewcommand{\thepage}{\arabic{page}}
\setcounter{page}{1}

\section{Introduction}%1

In a recent paper \cite{1} I showed how to expand 
\begin{equation}
F(x+\delta) = F(x) + F^{(1)}(x) + O(\delta^{2}) \label{a1}
\end{equation}
when both $x$ and $\delta$ were general quaternionic variables, thus 
did not commute with each other: 
\begin{equation}
F^{(1)}(x) = F^{\prime}(x) \delta_{1} + [F(x) - F(x^{*})] 
(x-x^{*})^{-1}\;\delta _{2}, \;\;\;\delta = \delta_{1} + 
\delta_{2},\label{a0}
\end{equation}
with specific formulas on how to construct the two components of 
$\delta$. 

Now we shall extend that analysis to a more general 
situation.

Consider the N x N matrices $X$ over the complex numbers and 
arbitrary analytic functions $F(X)$ with such a matrix as its 
variable. We seek 
a general construction for the first order term $F^{(1)}(X)$ when we 
expand $F(X + \Delta)$ given that $\Delta$ is small but still a general N x 
N matrix that does not commute with $X$.

The first step, as before, is to represent the function $F$ as a 
Fourier transform,
\begin{equation}
F(X) = \int dp f(p) e^{ipX}\label{a2}
\end{equation}
where the integral may go along any specified contour in the complex 
p-plane; and then we also make use of the expansion,
\begin{equation}
e^{(X+\Delta)} = e^{X}[1+\int_{0}^{1} ds e^{-sX}\Delta \;
e^{sX} + O(\Delta^{2})].\label{a3}
\end{equation}
Another well-known expansion, relevant to what we see in (\ref{a3}), 
is
\begin{equation}
e^{A}\;B\;e^{-A} = B + [A,B] + \frac{1}{2!} [A,[A,B]] + \frac{1}{3!}
[A,[A,[A,B]]] + \ldots  ,\label{a4}
\end{equation}
involving repeated use of the commutators, $[A,B]=AB-BA$.

\section{Diagonalization}%2

The first step is to assume that we can find a matrix $S$ that will 
diagonalize the matrix $X$ at any given point in the space of such 
matrices.
\begin{equation}
A = S\;X\;S^{-1}, \;\;\;\;\;\;\; A_{i,j} = \delta_{i,j}\; \lambda_{i}, 
\;\;\;\;\;\;\; i,j = 1,\ldots,N \label{b1}
\end{equation}
and we carry out the same transformation on the matrix $\Delta$:
\begin{equation}
B = S\;\Delta\;S^{-1}\label{b2}
\end{equation}
but, of course, the matrix $B$ will not be diagonal.

Our task is to separate the matrix $B$ into separate parts, each of 
which will behave simply in the expansion of Eq. (\ref{a4}).  The 
first step is to recognize that the diagonal part of the matrix B, 
call it $B_{0}$, commutes with the diagonal matrix $A$ and thus we 
have
\begin{equation}
e^{tA}\;B_{0}\;e^{-tA} = B_{0},\label{b3}
\end{equation}
where we have $t = -isp$ from (\ref{a2}), (\ref{a3}).

Next we look separately at each off-diagonal element of the matrix B: 
that is, $B_{(i,j)}$ is the matrix that has only that one 
off-diagonal ($i \neq j$) element of the given matrix $B$, and all the 
rest are zeros.  The first commutator is simply
\begin{equation}
[A,B_{(i,j)}] = r_{ij}\;B_{(i,j)}, \;\;\;\;\; r_{ij} = \lambda_{i} - 
\lambda_{j} \label{b4}
\end{equation}
and then the whole series can be summed:
\begin{equation}
e^{tA}\;B_{(i,j)}e^{-tA} = e^{tr_{ij}}\;B_{(i,j)}.\label{b5}
\end{equation}

Putting this all together, we have
\begin{equation}
S\;F^{(1)}\;S^{-1} = \int dp f(p)\;ip\;e^{ipA}\;\int_{0}^{1} ds\;[ 
B_{0} + \sum_{i \neq j}\;e^{-ispr_{ij}}\;B_{(i,j)}]. \label{b6}
\end{equation}

It is trivial to carry out the integrals over $s$; and we thus come 
to the final answer
\begin{equation}
 F^{(1)}(X) = F^{\prime}(X)\;\Delta_{0} + \sum_{i \neq j}
[F(X) - F(X-r_{ij}I)]\;r_{ij}^{-1}\;\Delta_{(i,j)} \label {b7}
\end{equation}
where
\begin{equation}
\Delta_{\mu} \equiv S^{-1}\;B_{\mu}\;S, \;\;\;\;\; \mu = 0, (i,j). 
\label{b8}
\end{equation}

\section{Discussion}%3

The general structure of the result, Eq.(\ref{b7}), is similar to 
what we found in earlier work, Eq.(\ref{a0}): the first term $F^{\prime}(X)$ looks 
like ordinary differential calculus and goes with that part of 
$\Delta$ that commutes with the local coordinate $X$; the remaining 
terms are non-local, involving the function $F$ evaluated at discrete 
points separated from $X$ by specific multiples of the unit matrix $I$.

While this final formula appears not as explicit as the previous 
result found for quaternionic variables (or for variables based upon 
the algebra of SU(2)), in any practical situation we have computer 
programs that can calculate the matrix operations referred to above 
with great efficiency.

The quantities $r_{ij}$, which may be real or complex numbers, can be 
called the ``roots'' following 
the familiar treatment of Lie algebras. They have some properties, such 
as $r_{ij} = - r_{ji}$ and sum rules that involve traces of the 
original matrix $X$ and powers of $X$.

What happens if the eigenvalues of $X$ are degenerate? Suppose, for 
example, that $\lambda_{1} = \lambda_{2}$. This means that $r_{12}$ 
and $r_{21}$ are zero. If we look at Eq. (\ref{b7}), we see that the 
terms $\Delta_{(1,2)}$ and $\Delta_{(2,1)}$ then have the coefficient 
$F^{\prime}(x)$.  Thus they simply add in with $\Delta_{0}$.  In the 
extreme case when all the eigenvalues are identical, then the answer 
is $F(X + \Delta) = F(X) + F^{\prime}(X)\;\Delta + O(\Delta^{2})$, 
which is old fashioned differential calculus for commuting variables.

\vskip 0.5cm
\setcounter{equation}{0}
\def\theequation{A.\arabic{equation}}
\boldmath
\noindent{\bf Appendix A}%A
\unboldmath
\vskip 0.5cm

A more abstract form of the result Eq. (\ref{b7}) is the following.
\begin{eqnarray}
&& F^{(1)}(X) = F^{\prime}(X)\;\Delta_{0}\; + [C,F(X)], \\ 
&& [C ,X] = \Delta - \Delta_{0}, \;\;\;\;\; [\Delta_{0},X] = 0. 
\end{eqnarray}
Here the matrix $C$ is defined implicitly, through its commutator 
with $X$, rather than explicitly; and 
the matrix $\Delta_{0}$ is the same as previously discussed. 
One may readily confirm the correctness of this formula in the case 
of $F(X) = X^{n}$.

This alternative formalism may also be applied to the case of a 
quaternionic variable $x$, which was studied in reference \cite{1}. In 
that case we find $\Delta_{0} = \delta_{1}$ and $C = \frac{1}{4r^{2}}
[x,\delta]$. 

It is interesting that this new formalism manages to hide the 
non-locality, which was a prominent feature of the original analysis.

\end{document}